\RequirePackage{ifpdf}
\ifpdf % We are running pdfTeX in pdf mode
\documentclass[pdftex]{sigma}
\else
\documentclass{sigma}
\fi

\newcommand{\ZZ}{\mathbb Z}    % Integer
\newcommand{\CC}{\mathbb C}    % Complex
\newcommand{\B}{\mathcal B}
\newcommand{\gl}{\mathfrak {gl}}    % Lie algebra gl

\begin{document}
\allowdisplaybreaks

\renewcommand{\thefootnote}{$\star$}

\renewcommand{\PaperNumber}{100}

\FirstPageHeading

\ShortArticleName{Sklyanin Determinant for Ref\/lection Algebra}

\ArticleName{Sklyanin Determinant for Ref\/lection Algebra\footnote{This paper is a
contribution to the Proceedings of the International Workshop ``Recent Advances in Quantum Integrable Systems''. The
full collection is available at
\href{http://www.emis.de/journals/SIGMA/RAQIS2010.html}{http://www.emis.de/journals/SIGMA/RAQIS2010.html}}}

\Author{Natasha ROZHKOVSKAYA}

\AuthorNameForHeading{N. Rozhkovskaya}

\Address{Department of Mathematics, Kansas State University, USA}
\Email{\href{mailto:rozhkovs@math.ksu.edu}{rozhkovs@math.ksu.edu}}

\ArticleDates{Received September 21, 2010, in f\/inal form December 23, 2010;  Published online December 29, 2010}

\Abstract{Ref\/lection algebras  is a class of algebras associated with integrable models with boundaries.  The coef\/f\/icients of Sklyanin determinant  generate the center of the ref\/lection algebra. We give a
combinatorial  description of Sklyanin determinant  suitable for explicit computations.}


\Classification{05E10; 17B37}

\renewcommand{\thefootnote}{\arabic{footnote}}
\setcounter{footnote}{0}

%\section{Introduction}
In \cite{Skl}  E.K.~Sklyanin introduced a class of algebras  associated with integrable models with bounda\-ries. Following \cite{MR}, we call them {\it reflection algebras}. The family  of ref\/lection algebras $\B(n,l)$ is def\/ined as associative algebras whose generators satisfy two types of  relations: the ref\/lection equation  and the unitary condition.
In \cite{MR}   A.I.~Molev and  E.~Ragoucy show that the center of $\B(n,l)$  is generated by the coef\/f\/icients  of an analogue of quantum determinant, called the {\it Sklyanin determinant}. In the same paper
 the authors develop an analogue of Drinfeld's highest  weight theory  for ref\/lection algebras and give a complete description of their f\/inite-dimensional  irreducible representations.

 The    ref\/lection algebras $\B(n,l)$  have many common features with the  {\it twisted Yangians}, introduced by G.~Olshanski \cite{Olsh}.
 For example, the center of twisted Yangian   is  also generated by coef\/f\/icients of Sklyanin determinant, which is  def\/ined similarly to the  Sklyanin determinant of ref\/lection algebras.
 The detailed exposition of the representation  theory of twisted Yangians can be found in~\cite{Mol1, Mol}.

 In \cite{Mol3} A.I.~Molev
gives a combinatorial formula  for the Sklyanin determinant of twisted Yangians in terms of matrix elements of the matrix of generators of the twisted Yangians (see also~\cite{Mol1}).
 The goal of this paper is to give  analogues combinatorial description of the Sklaynin determinant of the ref\/lection algebra.

 In Section~\ref{section1} we review the main def\/initions. In Sections~\ref{section2},~\ref{section3} we rewrite the def\/inition of the Sklyanin determinant of $\B(n,l)$
in an alternative form. We  describe combinatorics of the corresponding product of generating matrices `twisted' by Jucys--Murphy elements. In Section~\ref{section4}, Theorem~\ref{sdet}, we prove combinatorial formula for Sklyanin determinant in terms of matrix elements of the generating matrix of $\B(n,l)$.  In Section~\ref{section5} we give examples of computations by this formula.

\section{Def\/initions}\label{section1}
The following  notations will be used through the paper.
For a matrix  $X$ with entries $(x_{ij})_{i,j=1,\dots, n}$ in an  associative algebra $A$ write
\[
X_s=\sum_{i,j}1\otimes\dots\otimes \underset{s}{E_{ij}}\otimes\dots \otimes 1\otimes x_{ij} \in\operatorname{End}(\CC^n)^{\otimes k}\otimes A.
\]
 Consider the Yang matrix
\[
R(u)=1-\frac{P}{u}  \in \text{End}\,(\CC^n)^{\otimes 2}\big[\big[u^{-1}\big]\big],
\]
where $P:v\otimes w\to w\otimes v$ is a permutation operator in the space $\CC^n\otimes\CC^n$.

 For $i$, $j$ such that $1\le i < j\le k$ the notations $P_{ij}$ and
  $R_{ij}(u)$   are used  for the  action $P$ and  $R(u)$ on the $i$-th and $j$-th copies of  the vector space $\CC^n$ in $(\CC^n)^{\otimes k}$.

 Let $m$, $l$ be such nonnegative integers that $n=m+l$. Let $(\varepsilon _1,\dots, \varepsilon_n)$ be an $n$-tuple with entries $\pm 1$, such that
\[
 \varepsilon_i=
 \begin{cases}
\phantom{-}1& \text{if}\quad 1\le i\le m,\\
 {-}1& \text{if}\quad n\ge i> m.
 \end{cases}
\]

 The ref\/lection algebra $\B(n,l)$ is a unital associative algebra with the generators $b_{ij}^{(r)}$, $i,j=1,\dots, n$, $r= 0,1, 2,\dots $, that satisfy the  relations of the ref\/lection equation (\ref{RE}) and the unitary condition (\ref{un}) described below.

Combine the generators $b_{ij}^{(r)}$ into formal series
\[
b_i^{j}(u)=\sum_{r=0}^{\infty}b_{ij}^{(r)}u^{-r}\qquad \text{with }\quad
b^{(0)}_{ij}=\delta_{ij}\varepsilon_{i},
\]
and collect them into generating matrix
\[
B(u)=\sum_{i,j=1}^{n} E_{ij}\otimes b_i^{j}(u)\in \text{End}\, \CC^n\otimes \B(n,l)\big[\big[u^{-1}\big]\big].
\]
Then the ref\/lection equation  relation is given by
\begin{gather}\label{RE}
R(u-v)B_1(u)R(u+v)B_2(v)=B_2(v)R(u+v)B_1(u)R(u-v),
\end{gather}
and the unitary condition  is
\begin{gather}\label{un}
B(u)B(-u)=1.
\end{gather}

\begin{remark}  The unitary condition allows to realize the algebra $\B(n,l)$ as a subalgebra of Yangian  $Y(\gl_n(\CC))$. The mentioned below Proposition \ref{MR} is proved in \cite{MR}  through this inclusion of algebras.
However, we do not need the unitary conditions anywhere for the proofs of the statements below. Hence all the results of the paper automatically become true for the extended version of the ref\/lection algebra, which is def\/ined the same way as the algebra  $\B(n,l)$, but with the unitary condition (\ref{un}) omitted.
\end{remark}

We abbreviate  $R_{ij}:=R_{ij}(2u-i-j+2)$.
Put
\[
\langle B_1,\dots, B_k\rangle =B_1(u)(R_{12}\cdots R_{1k})B_2(u-1) (R_{23}\cdots R_{2k})
\cdots B_ k(u-k+1).
\]
This expression  is as an  element of $\text{End}\, (\CC^n)^{\otimes k }\otimes \B(n,l)[[u^{-1}]]$.

Let $A_n$ be the full  anti-symmetriztion operator in the space $\text{End}\, (\CC^n)^{\otimes n }$:
\begin{gather*}
A_n (v_1\otimes \dots\otimes v_n)=\sum_{\sigma\in S_n} (-1)^{\sigma}
v_{\sigma(1)}\otimes \dots \otimes v_{\sigma(n)}.
\end{gather*}
One can show \cite{MR} that $A_n\langle B_1,\dots, B_n\rangle $ is equal to the product of the anti-symmetrizer $A_n$ and a series in $u^{-1}$   with coef\/f\/icients in $\B(n,l)$. This series is called Sklyanin determinant.

\begin{definition} Sklyanin determinant of $B(u)$ is such a  series in $u^{-1}$   with coef\/f\/icients in $\B(n,l)$ that
\[
 A_n \langle B_1,\dots, B_n\rangle =\text{sdet} B(u)A_n.
 \]
 \end{definition}
\begin{proposition}[\cite {MR}]\label{MR} The coefficients of the series $\operatorname{sdet}B(u)$ are  central in $\B(n,l)$. The odd coefficients $c_1,c_3,\dots$ in the expansion
\[
\operatorname{sdet} B(u)=(-1)^l+c_1u^{-1}+c_2u^{-2}+c_3u^{-3}+\cdots
\]
 are algebraically independent generators of the center of  $\B(n,l)$.
 \end{proposition}

\section{Alternative def\/inition of  the Sklyanin determinant}\label{section2}

The following proposition allows to simplify the combinatorial description of $\text{sdet}\, B(u)$.
\begin{proposition} \label{new}
For $k=1,\dots, n,$ define
\[
\Pi_k:=1-{(2u-k-n+2)^{-1}}\sum_{i=k+1}^{n}P_{ki},
\]
and put
\begin{gather*}
\langle \langle B_1,\dots, B_n
\rangle \rangle :=
B_1(u) \Pi_1B_2(u-1)\Pi_2\cdots B_{n-1}(u-n+2)\Pi_{n-1}B_n(u-n+1).
\end{gather*}
Then
\begin{gather*}
A_n\langle \langle B_1,\dots, B_n
\rangle \rangle =\operatorname{sdet}B(u)A_n.
\end{gather*}
\end{proposition}
\begin{proof}
We use the notation  $A_k$  for the antisymmetrizator in  $(\CC^n)^{\otimes n}$ that acts on the f\/irst $k$ copies of $\CC^n$, and
$A^{\prime}_{n-k}$  for  the antisymmetrizator  in   $(\CC^n)^{\otimes n}$ that acts on  the last $(n-k)$ copies of $\CC^n$:
\begin{gather*}
A_k (v_1\otimes \cdots\otimes v_n)=\sum_{\sigma\in S_k} (-1)^{\sigma}
v_{\sigma(1)}\otimes \cdots \otimes  v_{\sigma(k) }\otimes v_{k+1}\otimes\cdots
\otimes v_{n},
\\
A^\prime_{n-k} (v_1\otimes \cdots\otimes v_n)=\sum_{\sigma\in S_{n-k}} (-1)^{\sigma}
v_{1}\otimes \cdots \otimes v_{k }\otimes v_{\sigma(k+1)}\otimes\cdots
\otimes v_{\sigma(n)}.
\end{gather*}

\begin{lemma} \label{lem_ARR}
\begin{gather}\label{ARR}
A^\prime_{n-k} R_{kk+1}\cdots R_{kn} =
A^\prime_{n-k} \Pi_k=
\Pi_k\,A^\prime_{n-k}.
\end{gather}
\end{lemma}
\begin{proof} We follow the lines of the proof of Proposition~1.6.2 in~\cite{Mol1}. Let  $k+1\le i_1< i_2< \cdots < i_m\le n$. Then
\[
(k,i_1)(k,i_2)\cdots (k,i_m)=(i_1,i_2)(i_1,i_3)\cdots (i_{1},i_m)(k,i_1),
\]
where $(a,b)$   denotes  a  transposition.
Therefore,
\begin{gather}\label{APP}
A^\prime_{n-k}P_{ki_1}P_{ki_2}\cdots P_{ki_m}=(-1)^{m-1}A^\prime_{n-k}P_{ki_1}.
\end{gather}
The product of $R$-matrices $R_{kk+1}\cdots R_{kn}$  is a sum of  the identity operator and of monomials $P_{ki_1}P_{ki_2}\cdots P_{ki_m}$ with rational coef\/f\/icients of the form
\[
\frac{(-1)^m}{(2u-i_1-k+2)\cdots (2u-i_m-k+2)}.
\]
Using (\ref{APP}), we substitute each monomial by  $P_{ki_1}$   and   collect the coef\/f\/icients  for each  of $i=k+1,\dots, n$. Then the  f\/irst equality of the lemma follows.

In the group algebra of the symmetric group  $S_n$ for any $i,j,k\in\{1,\dots, n\}$, such that $i$, $j$, $k$ are pairwise distinct numbers, one has
\begin{gather}\label{ijk}
(i,j)((k,i)+(k,j))=((k,i)+(k,j))(i,j).
\end{gather}
Since $A^{\prime}_{n-k}$ is a linear combination of products of  transpositions $P_{ij}$ with  $i>k$ and $ j>k $,  by (\ref{ijk})  the antisymmetrizer  $A^{\prime}_{n-k}$ commutes with the sum $\sum_{i=k+1}^{n}P_{ki}$, and the second equality follows.
\end{proof}

 For  $i\le m\le n$ the antisymmetrizer $A^{\prime}_{n-m}$ commutes with $B_i(u-i+1)$,
and by (\ref{ARR}) it also commutes with $\Pi_m$.
Using that for $m=2, 3,\dots, n $,
\[
(m-1)!A^{\prime}_m=A^{\prime}_mA^{\prime}_{m-1},
\]
and Lemma \ref{lem_ARR}, we can   rewrite the product  $A_n \langle  B_1, \dots, B_n\rangle $ in  the def\/inition of $\text{sdet}\, B(u)$ in terms of operators $\Pi_1,\dots, \Pi_n$:
 \begin{gather*}
 A_n \langle  B_1, \dots, B_n\rangle
=  { \prod_{i=1}^{n-1}\frac{1}{(n-i)!}}
   A_n\,B_1(u) \Pi_1A^{\prime}_{n-1}B_2(u-1)\Pi_2\dots A^{\prime}_{1}
B_n(u-n+1)\\
\phantom{A_n \langle  B_1, \dots, B_n\rangle   }{}
=  A_n B_1(u) \Pi_1B_2(u-1)\Pi_2 \cdots
B_n(u-n+1)=A_n \langle \langle  B_1, \dots, B_n\rangle \rangle .\tag*{\qed}
\end{gather*}
\renewcommand{\qed}{}
\end{proof}

\section[Properties of $\langle \langle  B_1, \dots, B_n\rangle \rangle $]{Properties of $\boldsymbol{\langle \langle  B_1, \dots, B_n\rangle \rangle }$}\label{section3}

Proposition~\ref{new}  states that Sklyanin determinant  can be computed  using the  antisymmetrized product  $\langle \langle  B_1, \dots, B_n\rangle \rangle $ instead of the antisymmetrized  product   $\langle  B_1, \dots, B_n\rangle $. In this section   we  give three  combinatorial  descriptions of the product $\langle \langle  B_1, \dots, B_n\rangle \rangle $.

Let $I_n\subset \ZZ^n$  be the subset of such $n$-tuples $\kappa=(k_1, \dots, k_n)$ that
\[
i\le k_i\le n \qquad \text{for all}\quad i=1,2,\dots, n.
\]
 For $\kappa\in I_n$ introduce  a rational function $\alpha{(\kappa)}=\alpha{(\kappa,u)}$ of   a variable $u$ def\/ined by
\begin{gather*}%\label {alpha}
\alpha{(\kappa)}=
{\prod_{i< k_{i}} \frac{1}{(n-2u-2+i)}}.
\end{gather*}

For all $k=1,\dots, n$ we identify $P_{kk}$ with $\text{Id}^{\otimes n}\in \text{End}\, (\CC^n)^{\otimes n}$. Then
\begin{gather}\label{bp}
\langle \langle
B_1,\dots, B_n
\rangle \rangle =
 \sum _{\kappa=(k_1,\dots, k_n) \in I_n} \alpha(\kappa)\,B_1(u) P_{1k_1}
\cdots
B_{n}(u-n+1) P_{nk_n}.
\end{gather}
(The last operator $P_{n k_n}=P_{nn}=\text{Id}^{\otimes n}$ in the end of the formula  is  added  for a uniform presentation.)

\begin{lemma}\label{PXXP} For any $X\in \operatorname{End}( \CC^n) \otimes\B(n,l)$
\[
P_{ij}X_i=X_jP_{ij},
\qquad
X_iP_{ij}=P_{ij}X_j.
\]
\end{lemma}
Using this simple fact   all permutation operators  in the expression
$\langle \langle B_1,\dots, B_n\rangle \rangle $  can be moved in front of all elements $B_s(u-a)$.

\begin{proposition}\label{QA}
 Put
 \[
 Q^{-1}_{\kappa}=P_{1 k_1}\cdots P_{n k_n},
 \]
and denote the corresponding element of symmetric group $S_n$ as
\[
q_{\kappa}=(n,k_n)(n-1,k_{n-1})\cdots (1,k_1)\in S_n.
\]
Define also
\[
q_{\kappa}^{[i]}=(n,k_n)(n-1,k_{n-1})\cdots (i,k_i)\in S_n.
\]
Then
\begin{gather}\label{q2}
\langle \langle B_1,\dots, B_n\rangle \rangle
=\sum_{ \kappa \in I_n }\alpha{( \kappa)}\,Q^{-1}_{\kappa}B_{q_{\kappa}^{[1]}(1)}(u)\cdots B_{q_{\kappa}^{[n]}(n)}(u-n+1).
\end{gather}
Here $q_{\kappa}^{[i]}(j)$ is the result of action of permutation $q_{\kappa}^{[i]}$ on $j$.
\end{proposition}

\begin{proof}
From Lemma \ref{PXXP}  and the equality (\ref{bp}) the statement of the proposition  follows imme\-dia\-te\-ly.
\end{proof}

For the second  combinatorial description of $\langle \langle B_1,\dots, B_n\rangle \rangle $ we assign to  each $\eta\in I_n$ a permutation
 $p_\eta$ that is def\/ined by the following recursive rule.

  Def\/ine a sequence of permutations $\{p^{[i]}_\eta, i=1,\dots, n\}$ by
  \begin{gather*}%\label{p}
    p^{[n]}_\eta :=\text{Id},\qquad
       p^{[i]}_\eta :=
       p^{[i+1]}_\eta \big(i,\big(p^{[i+1]}_\eta\big)^{-1}(\eta_i)\big)\qquad \text{for}\quad i=n-1, n-2,\dots, 1.
 \end{gather*}
Put $p_\eta:=p^{[1]}_\eta$.

\begin{proposition}\label{QB}
 Then
\begin{gather*}%\label{q1}
\langle \langle B_1,\dots, B_n\rangle \rangle
=\sum_{ \eta\in I_n }\alpha{( \eta)}\,p_\eta^{-1}B_{ \eta_1}(u)\cdots B_{ \eta_n}(u-n+1),
\end{gather*}
where $p_\eta^{-1}$   is an  element of the symmetric group $S_n$, identified with the corresponding operator acting on $(\CC^n)^{\otimes n}$.
 \end{proposition}

\begin{proof}
Consider $\kappa=(k_1, \dots, k_n)\in I_n$.  Put
\[
\eta=\eta(\kappa)=(\eta_1, \dots, \eta_n)
\] with
\begin{gather}\label{kn}
\eta_i:=q_{\kappa}^{[i]}(i)=(n,k_n)\cdots (i, k_i) (i)  \qquad\text{for}\quad  i=1,\dots, n,
\end{gather}
where $q_{\kappa}^{[i]}$ is from Proposition~\ref{QA}.
Then obviously $\eta_i\!\ge\! i$ for all $i=1, \dots, n$, so \mbox{$\eta\!=\!(\eta_1, \dots, \eta_n)\!\in\! I_n$}. Moreover, any $\eta\in I_n$ can be obtained from some $\kappa \in I_n$ by the formula~(\ref{kn}): given $\eta=(\eta_1, \dots, \eta_n)\in I_n$, def\/ine  such $\kappa=(k_1, \dots, k_n)$ recursively:
\begin{gather*}
k_n:= \eta_n,\qquad
k_i:= (i+1,k_{i+1})\cdots (n, k_n)\,(\eta_i)\qquad \text{for}\quad i=n-1,\dots, 2,1.
\end{gather*}
Therefore, the map $\kappa\mapsto \eta(\kappa)$ is a bijection from $I_n$ to~$I_n$. Set $p_{\eta(\kappa)}:=q_\kappa$, and  observe that $\alpha(\kappa)= \alpha (\eta(\kappa))$. Changing the index of summation $\kappa \in I_n$ to  $\eta \in I_n$ in  (\ref{q2}) proves the statement  of the Proposition~\ref{QB}.
\end{proof}

In Propositions \ref{QA},  \ref{QB} the sum of monomials  is taken over the set $I_n$. This set   has  $n!$ elements. One can ask, if the monomials can be  naturally numerated by permutations  $\sigma\in S_n$.

We will need the following notion of a {\it word restriction of a permutation}.
 Let $\sigma$  be a permutation of the elements of the set $\{1, \dots, n\}$, and let $\sigma=\Gamma_1\cdots \Gamma_l$ be a decomposition  into non-intersecting cycles $\Gamma_i=(g_{i,1},\dots, g_{i,s_i})$.  Let $G$ be a subset of $\{1, \dots, n\}$.

\begin{definition}
We say that the cycle $\Gamma |_G$ is {\it {the word restriction of  the cycle}}  $\Gamma$ on the set $G$, if~$\Gamma |_G$ is obtained from $\Gamma$ by deleting   from the  word of the cycle of $\Gamma$ of all the elements that do not belong to  $G$.
Then {\it {the word restriction}} of $\sigma$ on the set $G$ is def\/ined as a product of word restrictions of cycles $\Gamma_1|_G, \dots, \Gamma_l|_G $ on the set $G$:
\[
\sigma|_G=
\Gamma_1|_G \cdots \Gamma_l|_G.
\]
We say that a cycle $\Gamma$  has {\it{an empty restriction}} on the set $G$, if it contains no elements from $G$.
 \end{definition}
 \begin{example} The word restriction of the permutation $\sigma =(1,5,7,3)(4)(6,2)$ of the set $\{1,2,3,4,$ $5,6,7\}$ on the subset $G=\{1,2,3\}$
 is  the permutation $(1,3)(2)$.
 \end{example}

\begin{proposition}\label{QC}  For any  $\sigma\in S_n$ consider  a decomposition of $\sigma$  into nonintersecting cycles:
$\sigma = \Gamma_1\cdots \Gamma_s$.
Let $g_i$ be the maximal element of the cycle $\Gamma_i$.   Denote as $G_\sigma$ the set $\{g_i\}_{i=1}^{s}$.
Then
\[
\langle \langle B_1,\dots, B_n\rangle \rangle
=\sum_{\sigma \in S_n }\bar\alpha{(\sigma)} \sigma^{-1} B_{\sigma^{[1]}(1)}(u)\cdots B_{\sigma^{[n]}(n)}(u-n+1).
\]
Here
\begin{gather*}%\label {bar_alpha}
\bar\alpha{(\sigma)}=
\prod_{g\notin G_\sigma} \frac{1}{ (n-2u-2+g)},
\end{gather*}
the permutation $\sigma^{[i]}$ is the word restriction of  $\sigma$ on the set $\{i, \dots, n\}$. We identify $\sigma^{-1}$  with the corresponding operator acting on $(\CC^n)^{\otimes n}$.
\end{proposition}

\begin{proof}
Observe that any $\sigma\in S_n $ can be written  uniquely as a product of transpositions
\begin{gather}\label{gammak}
\sigma=(n,k_n)\cdots (1,k_1)
\end{gather}
 for some $\kappa=(k_1,\dots, k_n)\in I_n$. In other words, there  is a one-to-one correspondence between the elements of $I_n$  and permutations in $S_n$. One can check that the word restriction  of $\sigma$ on $\{i,\dots,n\}$  has the property
\[
\sigma^{[i]}=(n,k_n)\cdots (i,k_i)=q^{[i]}_{\kappa}.
\]
The set of   maximal elements  of the cycles  $\{g_{i}\}_{i=1}^{s}$  in the decomposition of $\sigma$ into non-intersecting cycles  coincides with the set $\{k_i\,| \, k_i=i\}$ in  the decomposition~(\ref{gammak}). This implies the formula for~$\bar\alpha(\sigma)$.
\end{proof}

 \section{Computation of Sklyanin determinant}\label{section4}
 In this  section we give formulas that allow to compute explicitly Sklyanin determinant in terms of matrix elements of the generating matrix $B(u)$.

For any  $n$-tuple $\eta\in I_n$  we will associate several quantities and sets. The   $n$-tuple  $\eta\in I_n$ can be  considered as a function on the interval of integer numbers $\{1,\dots,n\}$. If $\eta=(\eta_1,\dots, \eta_n)$, then
\[
\eta: \ \ \{1,\dots,n\}\to \{1,\dots, n\},\qquad \eta(i):=\eta_i.
\]
Let  the image $\text{Im} (\eta)$ consist of  distinct numbers $ \{N_1,\dots, N_K\}$,  and let the pre-image
$\eta^{-1}(N_l)=\{i\,|\,\eta_{i}=N_l \}$     consist of numbers  $\{a_{l,1}, \dots, a_{l,m_l}\}$. We choose the order of the elements so that
\[
a_{l,1}<\cdots<a_{l,m_l}.
\]
For each $l=1, \dots, K$ the sequence  $a_{l,1}<\cdots< a_{l,m_l}$ def\/ines a increasing cyclic permutation $\Gamma_l=(a_{l,1},\dots,a_{l,m_l}) \in S_n$. Put
\[
\gamma_\eta=\Gamma_1\cdots \Gamma_K.
\]
Thus, for every $\eta \in I_n$ we associate a permutation $\gamma_\eta\in S_n$ that is a product of  non-intersecting increasing cycles.
We also denote as $G_\eta^{\pm}$ the sets of maximal and minimal elements of those cycles:
\[
G_\eta^-=\{a_{l,1}\}_{l=1}^{K}, \qquad G_\eta^+=\{a_{l,m_l}\}_{l=1}^{K}.
\]

Let $S(\eta)$ be the the subgroup of  $S_n$  of all permutations that act trivially on the complement of the set $\operatorname{Im}(\eta)$ in  $\{1,\dots, n\}$.

 Consider the set $\{b^{p}_{q}(u)\}_{p, q=1,\dots, n}$ of matrix coef\/f\/icients of $B(u)$ and the  set of matrix elements  $\{b_{p_1, \dots , p_n}^{q_1, \dots ,q_n}(u,\eta)\}_{p_i, q_i = 1,\dots, n}$
 of $B_{ \eta_1}(u)\cdots B_{ \eta_n}(u-n+1)$:
\begin{gather*}
 B(u)=\sum_{p,q=1,\dots, n} E_{pq}\otimes b_p^q(u),
\\
B_{ \eta_1}(u)\cdots B_{ \eta_n}(u-n+1)=
\sum_{p_1,\ldots, p_n,\, q_1,\ldots, q_n}E_{p_1q_1}\otimes \cdots \otimes E_{p_nq_n}
b_{p_1, \dots, p_n}^{q_1, \dots, q_n}(u,\eta).
\end{gather*}

  \begin{theorem}\label{sdet}
\begin{gather*}%\label{sdet1}
\operatorname{sdet}\, B(u)=
\sum_{ \eta\in I_n } \sum_{\sigma\in S(\eta)}  \alpha(\sigma, \eta)
b^{1,\dots , n}_{\sigma(1),\dots, \sigma(n)}(u,\eta),
\end{gather*}
where $b^{1,\dots ,n}_{\sigma(1),\dots , \sigma(n)}(u,\eta)$ is the corresponding matrix coefficient of   $B_{\eta_1}(u) \cdots B_{\eta_n}(u-n+1)$. It can be computed by the following rule:
 \begin{gather}\label{b}
b^{1,\dots ,n}_{\sigma(1),\dots, \sigma(n)}(u,\eta)=
   \prod_{k=1}^{n}
b_{p(k)}^{q(k)}(u-k+1),\\
%\label{pq}
 p{(k)}=
\begin{cases}
\sigma(\eta(k))  &  \text{if}\quad k\in G_\eta^ -,\\
s_k& \text{otherwise},
\end{cases}
\qquad
q{(k)}=
\begin{cases}
\eta(k)  &\text{if}\quad
k\in G_\eta^+ ,\\
s_{\gamma_\eta(k)}&\text{otherwise},
\end{cases}\nonumber
\end{gather}
 and  the Einstein summation convention is used in \eqref{b} for the indices $s_t$: if  for some $t$ the repeated indices $s_t$ occur in \eqref{b}, the monomials are implicitly summed over  the range $s_t=1,\dots, n$.
Also
\[
\alpha(\sigma, \eta)=(-1)^{\sigma}
\prod_{i<\eta_i} \frac{1}{(2u+2-i-n)}.
\]
\end{theorem}

\begin{proof}
From Proposition \ref{QB}   we can write
\begin{gather}\label{eh1}
 \text{sdet}\, B(u) A_n=\left(
  \sum_{ \eta\in I_n } A_n
 \alpha{( \eta)}(-1)^{p_\eta}
  B_{ \eta_1}(u)\cdots B_{ \eta_n}(u-n+1)\right).
 \end{gather}

Following
the lines of the proof\footnote{The author would like to thank the referee for indicating this step.}
  of Proposition~2.7 in~\cite{MNO}
we apply both sides of (\ref{eh1}) to the vector $(e_1\otimes \cdots \otimes e_n)$. Let $v= A_n(e_1\otimes \cdots \otimes e_n)$. The left-hand side   of (\ref{eh1}) gives
$
 \text{sdet}\, B(u)v$.
The right-hand side equals
\[ \sum_{ \eta\in I_n } A_n
 \alpha{( \eta)}(-1)^{p_\eta}
\sum_{t_1, \dots ,  t_n}b^{1, \dots,  n}_{t_1,\dots, t_n}(u,\eta)
(e_{t_1}\otimes \cdots\otimes e_{t_n}).
\]
If not all of $\{t_1,\dots, t_n\}$ are pairwise distinct, then $A_n(e_{t_1}\otimes \cdots\otimes e_{t_n})=0$. Otherwise $e_{t_1}\otimes \cdots\otimes e_{t_n}=\sigma(e_1\otimes \cdots \otimes e_n)$, where the permutation $\sigma\in S_n$ is def\/ined by $\sigma(k)=t_k$, and
\[
A_n(e_{t_1}\otimes \cdots\otimes e_{t_n})=
A_n\sigma(e_1\otimes \cdots \otimes e_n)=(-1)^{\sigma}v.
\]

We obtain that
\begin{gather}
 \sum_{ \eta\in I_n }A_n
 \alpha{( \eta)}(-1)^{p_\eta}
\sum_{t_1,\dots ,  t_n}b^{1,\dots, n}_{t_1,\dots, t_n}(u,\eta)
(e_{t_1}\otimes \dots\otimes e_{t_n})\nonumber\\
\qquad =\sum_{\sigma\in S_n}\sum_{ \eta\in I_n }  \alpha{( \eta)}(-1)^{p_\eta}  (-1)^{\sigma}
b^{1,\dots, n}_{\sigma(1),\dots, \sigma(n)}(u,\eta)v.\label{rhs}
\end{gather}
Observe that if  for f\/ixed $\sigma \in S_n$ there exists $j\notin \text{Im}\,\eta$ such that $\sigma(j)\ne j$, then the corresponding matrix element $b^{1,\dots, n}_{\sigma(1),\dots, \sigma(n)}(u,\eta)$ vanishes to zero. Hence, we can substitute the sum in  (\ref{rhs}) over all permutations $\sigma \in S_n$ by the sum of permutations $\sigma\in S(\eta)$ and rewrite (\ref{rhs}) as
\begin{gather*}%\label{rhs1}
\sum_{ \eta\in I_n }  \sum_{\sigma\in S(\eta)}\alpha{( \eta)}(-1)^{p_\eta}  (-1)^{\sigma}
b^{1,\dots, n}_{\sigma(1),\dots, \sigma(n)}(u,\eta)v.
\end{gather*}
The direct check shows that $\alpha(\sigma, \eta)\!=\!\alpha{( \eta)}({-}\!1)^{p_\eta}\!  ({-}\!1)^{\sigma}$ and that the matrix coef\/f\/icients
$b^{1,{\dots}, n}_{\sigma(1),{\dots}, \sigma(n)}$ can be computed by~(\ref{b}).
\end{proof}

\section[Examples: $n=2,  3$]{Examples: $\boldsymbol{n=2,  3}$}\label{section5}

The Theorem \ref{sdet} allows to compute Sklyanin determinant explicitly.
\begin{example}
For $n=2$ the set $I_2$ has two elements: $I_2=\{(1,2), (2,2)\}$. Then
  $\eta=(1,2)$ contributes
  \[
  b_1^{1}(u)b_2^{2}(u-1) - b_2^{1}(u)b_1^{2}(u-1),
  \]
  and $\eta=(2,2)$ contributes
\[
\frac{1}{2u-1}\,(b_2^{1}(u)b_1^{2}(u-1)
+
b_2^{2}(u)b_2^{2}(u-1)),
\]
 so $\text{sdet}\,B(u)$ equals
\[
 b_1^{1}(u)b_2^{2}(u-1) - \frac{2u-2}{2u-1}b_2^{1}(u)b_1^{2}(u-1) + \frac{1}{2u-1}b_2^{2}(u)b_2^{2}(u-1).
 \]
 \end{example}

 \begin{example}
 For $n=3$ the set $I_3$ has six elements.
They make the following contributions to $\text{sdet} \,B(u)$:
\begin{gather*}
\eta=(1,2,3)\quad
 \sum_{\sigma\in S_3}(-1)^{\sigma}b^1_{\sigma(1)}(u)b^2_{\sigma(2)}(u-1)b^{3}_{\sigma(3)}(u-2),\\
\eta=(2,2,3)\quad
 \frac{1}{2u-2}\sum_{s=1}^3b^s_{2}(u)b^2_{s}(u-1)b^{3}_{3}(u-2)-
 b^s_{3}(u)b^2_{s}(u-1)b^{3}_{2}(u-2),\\
\eta=(3,2,3)\quad
\frac{1}{2u-2} \sum_{s=1}^3b^s_{3}(u)b^2_{2}(u-1)b^{3}_{s}(u-2)-
 b^s_{2}(u)b^2_{3}(u-1)b^{3}_{s}(u-2),\\
\eta=(1,3,3)\quad
 \frac{1}{2u-3}
 \sum_{s=1}^3b^1_{1}(u)b^s_{3}(u-1)b^{3}_{s}(u-2)-
 b^1_{3}(u)b^s_{1}(u-1)b^{3}_{s}(u-2),\\
\eta=(2,3,3)\quad
 \frac{1}{(2u-2)(2u-3)}
 \sum_{s=1}^{3}b^2_{2}(u)b^s_{3}(u-1)b^{3}_{s}(u-2)-
 b^2_{3}(u)b^s_{2}(u-1)b^{3}_{s}(u-2),\\
\eta=(3,3,3)\quad
 \frac{1}{(2u-2)(2u-3)}
  \sum_{s,t=1}^{3}b^s_{3}(u)b^t_{s}(u-1)b^{3}_{t}(u-2).
 \end{gather*}
The value  of $\text{sdet} \,B(u)$  equals the sum of these six  expressions.
We introduce the notation
\[
\left\vert\begin{matrix}
k&l&m\\
p&r&t
\end{matrix}\right\vert
:=b^k_{p}(u)b^l_{r}(u-1)b^{m}_{t}(u-2).
\]
Then
\begin{gather*}
\text{sdet} \,B(u)=
\left\vert\begin{matrix}
1&2&3\\
1&2&3
\end{matrix}\right\vert
+\frac{1}{2u-2}\left(
\left\vert\begin{matrix}
2&2&3\\
2&2&3
\end{matrix}\right\vert
+\left\vert\begin{matrix}
3&2&3\\
3&2&3
\end{matrix}\right\vert
-\left\vert\begin{matrix}
1&1&3\\
3&1&1
\end{matrix}\right\vert
-\left\vert\begin{matrix}
1&3&3\\
3&1&3
\end{matrix}\right\vert
\right)\\
\qquad{}+\frac{1}{2u-3}\left(
\left\vert\begin{matrix}
1&1&3\\
1&3&1
\end{matrix}\right\vert+
\left\vert\begin{matrix}
1&3&3\\
1&3&3
\end{matrix}\right\vert\right)
-\frac{2u-3}{2u-2}\left(
\left\vert\begin{matrix}
1&2&3\\
2&1&3
\end{matrix}\right\vert+
\left\vert\begin{matrix}
1&2&3\\
3&2&1
\end{matrix}\right\vert
-\left\vert\begin{matrix}
1&2&3\\
2&3&1
\end{matrix}\right\vert
\right)\\
\qquad {}-\frac{2u-4}{2u-3}
\left\vert\begin{matrix}
1&2&3\\
1&3&2
\end{matrix}\right\vert
-\frac{2u-4}{2u-2}
\left\vert\begin{matrix}
1&2&3\\
3&1&2
\end{matrix}\right\vert
+\frac{2u-4}{(2u-2)(2u-3)}\left(
\left\vert\begin{matrix}
2&2&3\\
2&3&2
\end{matrix}\right\vert+
\left\vert\begin{matrix}
3&2&3\\
3&3&2
\end{matrix}\right\vert
\right)
\\
\qquad{}+\frac{1}{(2u-2)(2u-3)}\left(
\left\vert\begin{matrix}
2&1&3\\
2&3&1
\end{matrix}\right\vert+
\left\vert\begin{matrix}
2&3&3\\
2&3&3
\end{matrix}\right\vert
+
\left\vert\begin{matrix}
3&1&3\\
3&3&1
\end{matrix}\right\vert
+
\left\vert\begin{matrix}
3&3&3\\
3&3&3
\end{matrix}\right\vert
\right).
\end{gather*}
 \end{example}
\vspace{-1mm}

\subsection*{Acknowledgements}
The author would like to thank A.I.~Molev and the referees   for valuable comments and suggestions that allowed to  simplify signif\/icantly the formulas and to improve the text of the paper. The travel to  RAQIS'10  was supported by AWM travel grant.
The research is partially supported by  KSU Mentoring fellowship for WMSE.

\vspace{-1mm}

\pdfbookmark[1]{References}{ref}
\LastPageEnding


\begin{thebibliography}{99}

\footnotesize\itemsep=-1pt

\bibitem{Mol1}
Molev A.,
Yangians and classical Lie algebras,
 {\it Mathematical Surveys and Monographs},  Vol.~143, American Mathematical Society, Providence, RI, 2007.

\bibitem{Mol}
Molev A.,
Yangians and their applications,
in  Handbook of Algebra,  Vol.~3, North-Holland, Amsterdam, 2003,  907--959,
\href{http://arxiv.org/abs/math.QA/0211288}{math.QA/0211288}.


 \bibitem{Mol3}
Molev A.,
Sklyanin determinant, Laplace operators, and charcteristic identities for classical Lie algebras,
\href{http://dx.doi.org/10.1063/1.531366}{{\it J.~Math. Phys.}} \textbf{36}  (1995), 923--943,
\href{http://arxiv.org/abs/hep-th/9409036}{hep-th/9409036}.

\bibitem{MNO}
Molev A., Nazarov M., Ol'shanskii G.,
Yangians and classical Lie algebras,
\href{http://dx.doi.org/10.1070/RM1996v051n02ABEH002772}{{\it Russian Math. Surveys}} \textbf{51}   (1996), no.~2, 205--282,
\href{http://arxiv.org/abs/hep-th/9409025}{hep-th/9409025}.

 \bibitem{MR}
 Molev A., Ragoucy E.,
Representations of ref\/lection algebras,
\href{http://dx.doi.org/10.1142/S0129055X02001156}{{\it Rev. Math. Phys.}} \textbf{14} (2002),  317--342,
\href{http://arxiv.org/abs/math.QA/0107213}{math.QA/0107213}.

 \bibitem{Olsh}
Ol'shanskii G.I.,
Twisted Yangians and inf\/inite-dimensional classical Lie algebras,
in  Quantum Groups (Leningrad, 1990),  Editor    P.~Kulish,  \href{http://dx.doi.org/10.1007/BFb0101183}{{\it Lecture Notes in Math.}}, Vol.~1510, Springer, Berlin, 1992, 104--119.


\bibitem{Skl}
Sklyanin E.K.,
Boundary conditions for integrable quantum systems,
\href{http://dx.doi.org/10.1088/0305-4470/21/10/015}{{\it J.~Phys.~A: Math. Gen.}} \textbf{21} (1988),  2375--2389.

\end{thebibliography}
 \end{document}